\documentclass[11pt]{amsart}
\usepackage[a4paper, total={6in, 8in}]{geometry}
\usepackage{epsfig}
\usepackage{amsmath}
\usepackage{amssymb}
\usepackage{bbold}
\PassOptionsToPackage{hyphens}{url}\usepackage{hyperref}

\usepackage{tikz}
\usetikzlibrary{plotmarks}
\usetikzlibrary{decorations.pathreplacing}
\usepackage{pgfplots}
\pgfplotsset{width=10cm,compat=1.9}
\usepgfplotslibrary{external}
\allowdisplaybreaks

\usepackage[labelfont=bf]{caption}
\usepackage{setspace}
\usepackage{subcaption}
\usepackage{cleveref}
\usepackage{lscape}
\usepackage{verbatim}
\usepackage{mathrsfs} 
\usepackage{multirow}
\usepackage{bm}
\usepackage{array} 
\usepackage{color}
\usepackage{bbm}
\newtheorem{theorem}{Theorem}[section]

\newtheorem{proposition}[theorem]{Proposition}

\usepackage{graphicx}
\usepackage{wrapfig}

\newcommand{\vb}{\vspace{3.2mm}}

\newcommand{\norm}[1]{\left\lVert#1\right\rVert}

\newcommand{\bs}{\boldsymbol}
\newcommand{\IA}{I_{\rm A}}
\newcommand{\IB}{I_{\rm B}}
\newcommand{\JA}{J_{\rm A}}
\newcommand{\JB}{J_{\rm B}}
\usepackage[normalem]{ulem}
\allowdisplaybreaks

\begin{document}

\title[A versatile stochastic dissemination model]{A versatile stochastic dissemination model}

\author[K.M.D. Chan, M.R.H. Mandjes]{K.M.D. Chan$^{1,2}$, M.R.H. Mandjes$^{1}$
}

\begin{abstract} 
This paper consider a highly general dissemination model that keeps track of the stochastic evolution of the distribution of wealth over a set of agents. There are two types of events: (i)~units of wealth externally arrive, and (ii)~units of wealth are redistributed among the agents, while throughout Markov modulation is allowed. We derive a system of coupled differential equations describing the joint transient distribution of the agents' wealth values, which translate into {\it linear} differential equations when considering the corresponding means and (co-)variances. While our model uses the (economic) terminology of wealth being distributed over agents, we illustrate through a series of examples that it can be used considerably more broadly. Indeed, it also facilitates the analysis of the spread of opinions over a population (thus generalizing existing opinion dynamics models), and the analysis of the dynamics of a file storage system (thus allowing the assessment of the efficacy of storage policies).

\vb


\vb

\noindent
{\sc Acknowledgments.} Dieter Fiems and Koen de Turck (Ghent University) are thanked for useful discussions.
The work in this paper was supported by the Netherlands Organisation for Scientific Research (NWO) through Gravitation-grant NETWORKS-024.002.003.

\vb

\noindent
$^{1}$Korteweg-de Vries Institute, University of Amsterdam, the Netherlands
    
\noindent  $^{2}$Transtrend {\sc bv}, Rotterdam, the Netherlands

\noindent    \vspace{0.1cm}
    
\noindent Corresponding author: Kit Ming Danny Chan ({\tt k.m.d.chan@uva.nl}).

\end{abstract}

\maketitle

\newpage 

\section{Introduction}
A stylized yet highly general model that describes the spread of wealth over a population of agents could encompass the following elements. In the first place, the agents are endowed with external inflow of wealth, for instance representing their salaries, paid to them by agents outside the population. In the second place, there are transactions: agents may purchase commodities or services from other agents, thus redistributing the wealth. Thirdly, to make the model more realistic, one could impose Markov modulation on the system: all parameters involved are affected by an exogenously evolving Markovian process. This Markovian background process could for instance represent the state of the economy, alternating between economic growth and recession.

The main contributions of this paper are the following. We develop a highly general dissemination model that keeps track of the stochastic evolution of the distribution of wealth over a set of agents, incorporating the three elements mentioned above (i.e., external inflow, redistribution and Markov modulation). A main asset of the model lies in it being broadly applicable and allowing closed-form analysis at the same time. 

For the resulting model, we succeed in deriving a system of coupled differential equations that describe the joint transient probability generating function of the agents' wealth levels, jointly with the state of the Markovian background process. When focusing on the corresponding means and (co-)variances, this system takes on a more convenient form, in that it becomes a system of {\it linear} differential equations, thus allowing for straightforward numerical evaluation. In passing, we also consider the model's stationary behavior, in particular establishing a stability condition. 

The broad applicability of the model is illustrated through a series of examples. As suggested by our terminology, being in terms of wealth that is distributed over agents, it can be used to analyze the evolution (in time) of a wealth vector. 
Another example concerns the dissemination of opinions over a population, where all agents influence one another. Our modelling framework extends existing opinion dynamics models in the way we incorporate stochasticity, while also the Markovian background process is a novel element. In a last example we consider a file storage system, where files of clients are periodically copied to one or multiple central storage units. Our model can be used to assess the efficacy of policies, intended to strike a proper balance between storage cost on one hand and the risk of data loss on the other hand. 

\vb

The dissemination model analyzed in this paper can be seen as a next step in a long tradition of queueing and population models.
In most existing queueing network models, the dynamics are such that the number of clients per queue changes by one at a time; see e.g.\ the accounts in \cite{KEL,SER}. A relatively small branch of the queueing literature considers queues with batch arrivals and batch services. In this respect we refer to e.g.\ 
\cite{COY,HEN,MIT} where product-form results are obtained. The setup that is probably closest to the one we consider in this paper, is the one of \cite{FMP18}, where the transaction events correspond to  the network population vector undergoing a (deterministic) linear transformation. In population-process theory and epidemics there is a strong emphasis on deterministic models to describe the dynamics of the sizes of various subpopulations (e.g.\ age groups, infected individuals, etc.); for more background, see for instance the monograph \cite{REN}. These models' stochastic counterparts have been studied as well; a general framework has been presented, and analyzed, in \cite{KUR}.
Importantly, to the best our knowledge, none of the earlier works covers our, highly general, redistribution mechanism.

This paper is organized as follows. The model and notation are introduced in Section~\ref{Sec:NOT}. Then subsequently the transient joint probability generating function (Section~\ref{Sec:PGF}), the first moments (Section~\ref{sec:FM}), and second moments ---~and hence also variances and covariances~---~(Section~\ref{sec:SM}) are analyzed; the section on first moments in addition establishes the model's stability criterion. Then there are three sections with illustrative examples, focusing on wealth redistribution (Section~\ref{Sec:WD}), on opinion dynamics (Section~\ref{sec:OD}) and on storage sharing systems (Section~\ref{Sec:SS}). Section~\ref{Sec:CR} concludes.

\section{Model and notation}\label{Sec:NOT}
In our model we study the stochastic behavior of ${\boldsymbol M}(t)\equiv (M_1(t),\ldots.M_I(t))$, where $M_i(t)$ denotes the ``wealth'' of agent $i$ at time $t$, with $i = 1, \ldots, I$, for some $I\in{\mathbb N}$. We recall that, as pointed out in the introduction,  ``wealth'' is to be interpreted in the broad sense; as we will extensively argue, the setup considered can also be used e.g.\ in the context of opinion spreading dynamics, or the context of file storage systems. 

To make our model as rich as possible, we let its dynamics be affected by an autonomously evolving Markovian background  {(or regime-switching)} process. In the (economic) context of wealth being spread over a population of individuals, the background process could reflect the state of the economy (e.g.\ alternating between economic peaks and periods of recession). Let this regime-switching process be modelled by the continuous-time Markov process $(X(t))_{t\geqslant 0}$  on the state space $\{1, \ldots, d\}$, for some $d\in{\mathbb N}$. This process, which is assumed to be irreducible, is governed by the transition rate matrix $Q = \{q_{ij}\}^{d}_{i,j=1}$ (with all non-diagonal elements being non-negative and row sums equal to 0), so that \cite{ASM03,NOR97}
\[{\mathbb P}(X(t) =\ell\,|\,X(0)=k) = \big(e^{Qt}\big)_{k,\ell}.\]

We proceed by describing the dynamics of the wealth process, given the background process is in state $k\in\{1,\ldots,d\}.$ We distinguish two types of events.
\begin{itemize}
    \item[$\circ$] In the first there are ``external arrivals'' of wealth. Concretely, for $j\in\{1,\ldots,J\}$ with $J\in{\mathbb N}$, at Poisson epochs with rate $\lambda_{jk}>0$ these external arrivals occur, leading to an increase of the wealth of all agents $i\in{S_j}\subseteq \{1,\ldots,I\}$ by one unit.   
    \item[$\circ$] In the second place there are ``shocks'', arriving to the system  according to a Poisson process with rate $\gamma_k>0$. At such a shock,  ``transactions'' take place, which concretely means that each of the $M_i(t)$ wealth units of agent $i$ contributes $W_{ijk} \in \mathbbm{N}_0$ wealth units to agent $j$.
    
    \noindent The precise mechanism is described more formally as follows. Supposing a shock happens at time $t>0$, then the number of wealth units after the shock at agent $j$ is, conditional on ${\bs M}(t-)  = (m_1, \ldots, m_I)^\top$ being the wealth vector just prior to time $t$, given by
    \[\sum_{i=1}^I \sum_{n=1}^{m_i} W_{ijkn},\]
    with $(W_{ijkn})_{n\in{\mathbb N}}$ denoting a sequence of independent and identically distributed random variables, all of them distributed as the discrete, non-negative random variable $W_{ijk}$.
    The random variables $W_{ijk}$ (with $i=1,\ldots I$ and $k=1,\ldots,d$)  are  assumed independent; {importantly, throughout we {\it do} allow dependence in $j$.}
     We define, for ${\bs z}$ such that $\max\{|z_1|,\ldots,|z_I|\}\leqslant {1}$, the associated probability generating function ({\sc pgf}) by
    \[ g_{ik}({\bs z}) = {\mathbb E}\Big(\prod_{j=1}^I z_j^{W_{ijk}} \Big).\]
\end{itemize}
Our primary aim is to establish a unique characterization of the transient wealth vector ${\bs M}(t)$, jointly with the state of the background process $X(t)$. To this end we work with the corresponding multivariate time-dependent joint {\sc pgf}.
Concretely, our analysis aims at identifying the following key object of study, for ${\bs z}$ such that $\max\{|z_1|,\ldots,|z_I|\}\leqslant {1}$:
\begin{equation} \label{fkzt}f_k(\boldsymbol{z},t) := {\mathbb E}\Big[\prod_{j=1}^{I}z_j^{M_j(t)} \mathbbm{1}_{\{X(t)=k\}}\Big],\end{equation}
which uniquely defines the distribution of $({\boldsymbol M}(t),X(t))\in {\mathbb N}^I\times\{1,\ldots,d\}.$ 

\section{Derivation of the differential equation for the joint PGF}\label{Sec:PGF}
The main objective of this section is to establish  a system of coupled differential equations (in $t$) for the time-dependent joint {\sc pgf}s $f_k(\boldsymbol{z},t)$, as defined in~\eqref{fkzt}. 
We do so relying on a standard argumentation: we  relate their values at time $t+\Delta t$ to their values at time $t$, with the aim to set up a system of differential equations. To this end, the underlying idea is to distinguish the three types of events that can occur in an interval of length $\Delta t$: a transition of the background process, external arrivals, and shocks (and, evidently, there is in addition the event that none of these three types of events occurs). Following this line of reasoning, we obtain for the time-dependent joint {\sc pgf} at time $t+\Delta t$ that
\begin{align*}
    f_k(\boldsymbol{z},t&+\Delta t) = \sum_{\ell\not = k}^dq_{\ell k}\,\Delta t \,{\mathbb E}\Big[\prod_{j=1}^{I}z_j^{M_j(t)} \mathbbm{1}_{\{X(t)=\ell\}}\Big]\:+\\
    \sum_{j=1}^J&\lambda_{jk}\,\Delta t \left(\prod_{i\in S_j} z_i\right) {\mathbb E}\Big[\prod_{i=1}^{I}z_i^{M_i(t)} \mathbbm{1}_{\{X(t)=k\}}\Big]\:+\\
    \gamma_k\,&\Delta t\sum_{{\boldsymbol m}\in{\mathbb N}^I} {\mathbb E}\Big[\prod_{j=1}^{I}z_j^{M_j(t+\Delta t)} \mathbbm{1}_{\{X(t+\Delta t)=k\}}\,\Big|\,{\boldsymbol M}(t)= {\boldsymbol m},{\mathscr E}_k(t)\Big]\,{\mathbb P}({\boldsymbol M}(t)= {\boldsymbol m})\:+\\
    \Big(1&-\sum_{\ell\not = k}^dq_{\ell k}\Delta t-\sum_{j=1}^J\lambda_{jk}\Delta t-\gamma_k\Delta t\Big) f_k(\boldsymbol{z},t)+o(\Delta t),
\end{align*}
with 
${\mathscr  E}_{k}(t)$ denoting the event of a shock between times $t$ and $t+\Delta t$ (evidently, while the background state is $k$). The right-hand side of the previous display can be interpreted and rewritten as follows. 
\begin{itemize}
    \item[$\circ$] The first term, which considers the scenario that the background process was in a state $\ell\not = k$ at time $t$ and makes a transition to $k$ between $t$ and $\Delta t$, equals by definition
    \[\sum_{\ell\not = k}^dq_{\ell k}\,\Delta t \,f_\ell(\boldsymbol{z},t).\]
    \item[$\circ$]The second terms represents the contributions of the external arrivals: if it is of the $j$-th type, then it increases (by 1) the wealth values of all agents $i$ such that $i\in S_j$. It reads
    \[\sum_{j=1}^J\lambda_{jk}\,\Delta t \left(\prod_{i\in S_j} z_i\right) f_k(\boldsymbol{z},t).\]
    \item[$\circ$]The third term describes the effect of the shocks. The claim is that we can express it in terms of the {\sc pgf} $f_k(\boldsymbol{z},t)$, but not evaluated in ${\bs z}$ but rather in a different argument. Indeed, again up to $\Delta t$-terms, due to the shock that occurs between times $t$ and $t+\Delta t$,
\begin{align*}
   \sum_{{\boldsymbol m}\in{\mathbb N}^I}{\mathbb E}&\Big[\prod_{j=1}^{I}z_j^{M_j(t+\Delta t)} \mathbbm{1}_{\{X(t+\Delta t)=k\}}\,\Big|\,{\boldsymbol M}(t)= {\boldsymbol m},{\mathscr E}_{k}(t)\Big]\,{\mathbb P}({\boldsymbol M}(t)= {\boldsymbol m})
   \\    &=\sum_{{\boldsymbol m}\in{\mathbb N}^I} \prod_{j =1}^I {\Big(g_{jk}({\bs z})\Big)^{m_j}}
   \,{\mathbb P}({\boldsymbol M}(t)= {\boldsymbol m},X(t)=k)\\
   &={\mathbb E}\Big[\prod_{j=1}^{I}g_{jk}({\bs z})^{M_j(t)} \mathbbm{1}_{\{X(t)=k\}}\Big]
   =
    f_k({\bs h}_k({\bs z}),t), 
\end{align*}
where we define
\[{\bs h}_k({\bs z}) := \big(g_{1k}({\bs z}),\ldots,g_{Ik}({\bs z})\big).\]
    \item[$\circ$]The fourth term corresponds to the scenario of no transition of the background process, no external arrivals and no shocks, leaving the wealth process unchanged. 
\end{itemize}

Observe that we have succeeded in expressing $f_k(\boldsymbol{z},t+\Delta t)$ in terms of quantities of the same type, as well as quantities of the type $f_\ell({\bs \varphi}({\bs z}),t)$ for known functions ${\bs \varphi}({\bs z}):[-1,1]^I\to [-1,1]^I.$
The next step is to subtract $f_k(\boldsymbol{z},t)$ from both sides of the equation, divide by $\Delta t$, and let $\Delta t\downarrow 0$, so as to obtain a system of differential equations in $t$. Indeed, we obtain, for $t\geqslant 0$,
\begin{equation} \label{eq_PGFDifToTime}
    \begin{split}
        \frac{\partial}{\partial t} f_k(\boldsymbol{z},t) =  \sum_{\ell=1}^d &q_{\ell k}f_\ell(\boldsymbol{z},t) \:+\\ \sum_{j=1}^J &\lambda_{jk}\left(\prod_{i\in S_j}z_i - 1\right) f_k(\boldsymbol{z},t) + \gamma_k \big(f_k({\bs h}_{k}(\boldsymbol{z}),t) - f_k(\boldsymbol{z},t) \big),
    \end{split}
\end{equation}
where we have used that the row sums of $Q$ are equal to $0$. Using the obvious property that ${\bs h}_k({\bs 1}) = {\bs 1}$, it is readily seen that 
$f_k(\boldsymbol{1},t)= {\mathbb P}\big(X(t) = k \big)$, as it should.
We have thus arrived at the following result.
\begin{proposition}\label{PRO:FDE}
    For any $t\geqslant 0$, $f_k(\boldsymbol{z},t)$ satisfies the system of differential equations \eqref{eq_PGFDifToTime}. If ${\bs M}(0) = {\bs m}_0$ and $X(0)=k_0$, then the initial condition is 
    \[f_k({\bs z},0)= \mathbbm{1}_{\{k=k_0\}} \prod_{i=1}^I z_i^{m_{0,i}}.\]
\end{proposition}
The above system of differential equations can be solved numerically. Importantly, the system is {\it not} linear; observe that in one term on the right-hand side of \eqref{eq_PGFDifToTime} the {\sc pgf} has the argument ${\bs h}_{k}(\boldsymbol{z})$ rather than ${\bs z}$. As we will notice in the next sections, however, when considering the computation of time-dependent moments (rather than the full time-dependent {\sc pgf}) we {\it do} obtain a reduction to systems of differential equations that {\it are} linear. These  can be solved at relatively low numerical effort using standard computational software. In the next two sections we subsequently concentrate on the evaluation of the first and second moments.

\section{Derivation of first moments, stability}\label{sec:FM}

In principle all moments of the components of ${\bs M}(t)$, as well as all mixed moments, can be derived from the differential equations \eqref{eq_PGFDifToTime} by repeated differentiation and plugging in $\boldsymbol{z}=\boldsymbol{1}$. In this section we discuss this widely applied procedure to determine the time-dependent first moments. Knowledge of these first moments also provides us with a criterion under which the model has a stable stationary version (i.e., does not explode as $t\to\infty$). 

\subsection{Differential equations for first moments} In this subsection we derive a system of linear differential equations that characterize the expectation of ${\bs M}(t)$.
With $w_{ijk}:={\mathbb E}W_{ijk}$ and $m_{ik}(t) := {\mathbb E}[M_i(t) \mathbbm{1}_{\{X(t)=k\}}]$, differentiating \eqref{eq_PGFDifToTime} to $z_i$  and inserting $\boldsymbol{z}=\boldsymbol{1}$ yields the following system of coupled ordinary differential equations:
\[m'_{ik}(t)= \sum_{\ell=1}^d q_{\ell k} m_{i\ell}(t) +
\sum_{j:i\in S_j}\lambda_{jk} \,\pi_k(t)+\gamma_k \left(\sum_{j=1}^I w_{jik}\, m_{jk}(t) - m_{ik}(t)\right)
,\]
with $\pi_k(t)={\mathbb P}(X(t)=k).$
Here we have used the standard differentiation rule for compositions of functions with vector-valued arguments, i.e., 
\[\frac{\partial f_k({\bs h}_{k}(\boldsymbol{z}),t)}{\partial z_i}=
\sum_{j=1}^I \frac{\partial f_k({\bs x},t) }{\partial x_j}\Big|_{{\bs x}={\bs h}_k({\bs z})} \frac{\partial (h_k({\bs z}))_j}{\partial z_i},\]
and
\[\frac{\partial (h_k({\bs z}))_j}{\partial z_i}\Big|_{{\bs z}={\bs 1}}=w_{jik}.\]
The next step is to compactly write the above system of differential equations  in matrix-vector form. 
We let ${\bs m}(t)\in{\mathbb R}^{dI}$ denote the stacked vector that results from
the $I$ vectors ${\bs m}_i(t)\equiv (m_{i1}(t),\ldots, m_{id}(t))^{\top}$, where ${\bs \pi}(t)\in{\mathbb R}^{dI}$ is  the stacked vector that results from
the $I$ (identical) vectors ${\bs \pi}_i(t)\equiv (\pi_{1}(t),\ldots, \pi_{d}(t))^{\top}$, i.e.,
\[{\bs m}(t) :=\left(\begin{array}{c}{\bs m}_1(t)\\
\vdots\\
{\bs m}_I(t)\end{array}
\right),\:\:\:\:{\bs \pi}(t) :=\left(\begin{array}{c}{\bs \pi}_1(t)\\
\vdots\\
{\bs \pi}_I(t)\end{array}
\right).\]
In addition, 
$G_{ji}:={\rm diag}\{\gamma_1 w_{ji1},\ldots,\gamma_d w_{jid}\}-{\rm diag}\{\gamma_1,\ldots,\gamma_d\}\mathbbm{1}_{\{i=j\}}$. We then define
\[A:=\left(\begin{array}{ccccc}
Q^\top+G_{11}&G_{21}&G_{31}&\cdots&G_{I1}\\
G_{12}&Q^\top+G_{22}&G_{{32}}&\cdots&G_{I2}\\
G_{13}&G_{23}&Q^\top+G_{33}&\cdots&G_{I3}\\
\vdots&\vdots&\vdots&\ddots&\vdots\\
G_{1I}&G_{2I}&G_{3I}&\cdots&Q^\top+G_{II}
\end{array}\right),\]
and, with $\Lambda_i:={\rm diag}\{\bar\lambda_{i1},\ldots,\bar\lambda_{id}\}$ and $\bar\lambda_{ik}:=\sum_{j:i\in S_j}\lambda_{jk}$,
\[\Lambda:=\left(\begin{array}{cccc}
\Lambda_1&0&\cdots&0\\
0&\Lambda_2&\cdots&0\\
\vdots&\vdots&\ddots&\vdots\\
0&0&\cdots&\Lambda_I\end{array}\right).\]
We thus end up with a system of $dI$ coupled linear differential equations, as given in the following proposition.
\begin{proposition} \label{PRO:MDE} For any $t\geqslant 0$, ${\bs m}(t)$ satisfies the system of linear differential equations
\[{\bs m}'(t) = A\,{\bs m}(t)+\Lambda\,{\bs \pi}(t).\]    
If ${\bs M}(0) = {\bs m}_0$ and $X(0)=k_0$, then the initial condition is 
    \[m_k(0)= \mathbbm{1}_{\{k=k_0\}} \,{m_{0,k}}.\]
\end{proposition}
This non-homogeneous system of linear differential equations can be solved in the standard manner.
In the first place, the vector ${\bs \pi}(t)$, corresponding to the transient state probabilities of the background process $X(t)$, satisfies the differential equation 
\[{\bs \pi}'(t) = ({\mathbb I} \otimes Q^{\top}) {\bs \pi}(t),\]
with $\otimes$ being the usual notation for the Kronecker product and ${\mathbb I}$ an identity matrix of appropriate dimension. This means, with
$\bar Q$ denoting the $(dI\times dI$)-dimensional matrix ${\mathbb I} \otimes Q^{\top}$, that ${\bs \pi}(t) = e^{\bar Q t}{\bs \pi}(0)$, or equivalently, 
\[\pi_j(t) = \sum_{i=1}^d {\mathbb P}(X(0)=i) \big( e^{Qt}\big)_{i,j}.\]
In the second place, the solution for ${\bs m}(t)$ can be written in terms of matrix exponentials, as follows:
\[{\bs m}(t) = e^{At}\,{\bs m}(0) + \int_0^t e^{A(t-s)} \,\Lambda\,{\bs \pi}(s)\,{\rm d}s.\]
We thus obtain the following result.
\begin{proposition} For any $t\geqslant 0$,
\[{\bs m}(t) = e^{At}\,{\bs m}(0) + \int_0^t e^{A(t-s)} \,\Lambda\,e^{\bar Q s}{\bs \pi}(0)\,{\rm d}s.\]
\end{proposition}

\subsection{Stability condition}
Let ${\bs \pi}:={\bs \pi}(\infty)$ be the unique solution of ${\bs \pi}Q={\bs 0}$ such that its entries sum to $1$, i.e., the stationary distribution of $X(t)$. 
In case the underlying model is stable, the above results directly imply that the steady-state mean vector ${\bs m}$ can be written in terms of the steady-state probabilities ${\bs \pi}$, as follows:
\[{\bs m} = - A^{-1} \Lambda {\bs \pi}.\]
The formal stability condition is given in the next statement. We define by $\omega$ the eigenvalue of $A$ with largest real part, i.e., the spectral abscissa of $A$. 
\begin{proposition}\label{pr:stab}
    The Markov chain $({\bs M}(t),X(t))_{t\geqslant 0}$ is ergodic if $\omega <0.$
\end{proposition}

\noindent {\it Proof.} The proof step-by-step mimics the one of \cite[Prop.\ 3]{FMP18}; we therefore restrict ourselves to sketching its main steps. The underlying idea is to establish  ergodicity of the skeleton Markov chain $({\bs M}(n\Delta),X(n\Delta))_{n\in{\mathbb N}}$ if $\omega <0$, for some $\Delta > 0$, where it should be noted that if the skeleton Markov chain is ergodic for some $\Delta > 0$, then so is $({\bs M}(t),X(t))_{t\geqslant 0}$ (observing that the mean recurrence time for any state of the skeleton chain is an upper bound for the mean
recurrence time of the original process). Then, by \cite[Prop.\ I.5.3]{ASM03}, a sufficient condition for ergodicity can be phrased in terms of 
\[{\mathbb E}\big[\norm{{\bs M}(\Delta)}_{1}\,|\,{\bs M}(0)={\bs m}_0, X(0)=k_0\big] -\norm{{\bs m}_0}_1<
 -\varepsilon,\]
 for some $\varepsilon>0$, and all ${\bs m}_0\in{\mathbb N}^I$ and $k_0\in\{1,\ldots,d\}$. Informally, this criterion entails that the process' drift is negative, bounded away from zero; cf.\ Foster's criterion \cite{BRE99,FOS53}. Such a bound can be achieved under $\omega<0$, with the precise same argumentation as the one used in the proof of \cite[Prop.\ 3]{FMP18}, where we rely on  \cite[Prop.\ 11.18]{BER09}  to find the required bound on the norm of
the matrix exponential. \hfill$\Box$
 
 \vb
 
 In case the stability condition $\omega < 0$ is not fulfilled, the components of ${\bs m}(t)$ typically grow in a fixed proportion; we return to this issue in Section \ref{sec:OD}. 
 
 \subsection{Special cases} In this subsection we provide more explicit results for three special cases: (1)~a single fully homogeneous population, (2)~a single distinct agent (a ``leader'') with $I-1$  homogeneous other agents (``followers''), (3)~two, internally homogeneous, interacting subpopulations. 
 
 \subsubsection{Homogeneous population.} We denote by $m_k(t)$ the mean wealth of an arbitrary agent, jointly with the event that the background process is in state $k$. As  a consequence of the fact that we ``start symmetrically'', i.e.,  the initial wealth of all agents is the same,  at any point in time the mean wealth of the individual agents coincides. In case we do not ``start symmetrically'' (for instance with two possible values of the initial wealth), the computations can still be performed, but become less clean. 
 
  In the variant that we consider, we let $J=I$, and we take $S_j=\{j\}$  and $\lambda_{jk}=\lambda_k$, for $j=1,\ldots,I$. (Here we note  that one can construct other fully symmetric external arrival processes, for instance by letting $J=1$ and  $S_1=\{1,\ldots,I\}$. Such alternative symmetric variants can be dealt with analogously.) The generic random variable $W_{ijk}$ now depends on the background state $k$ only (i.e., not on the indices $i$ and $j$ that indicate the agents); we let $w_k$ be its expected value.

 It is directly verified that our earlier results now yield
 \[m_k'(t) = \sum_{\ell=1}^d q_{\ell k}m_\ell(t) + \lambda_k \pi_k(t) +\gamma_k(Iw_k-1) \,m_k(t),\]
 or, in matrix-vector notation,
 \[{\bs m}'(t) = A\,{\bs m}(t)  + \Lambda \,{\bs \pi(t)},\]
 where $A:=Q^\top+{\rm diag}\{\gamma_1(Iw_1-1),\ldots,\gamma_d(Iw_d-1)\}$ and $\Lambda:={\rm diag}\{\lambda_1,\ldots,\lambda_d\}.$ {The $k$th entry of ${\bs m}$ now expresses the agents' mean wealth when the system is in state $k$.} In steady state, we obtain that the mean wealth vector equals $-A^{-1}\Lambda{\bs \pi}$, with ${\bs \pi}$ as defined before, provided that the stability condition is fulfilled (i.e., that the spectral abscissa of $A$ is negative). 
 
  \subsubsection{Leader and homogeneous followers.} \label{ssec:LHF} In this model, there is a single leader
 and $I-1$ homogeneous followers. 
 We ``start symmetrically'', i.e., all followers have the same initial wealth. We let $m_{{\rm L},k}(t)$ be the mean wealth of the leader at time $t$, and $m_{{\rm F},k}(t)$ the mean wealth of an arbitrary follower at time $t$, both jointly with the event that the background process is in state $k$.
 
As we did in the case of a homogeneous population, we take $J=I$ with $S_j=\{j\}$. We let the external arrival rate of the leader be $\lambda_{{\rm L},k}$, and of the followers $\lambda_{{\rm F},k}$, both corresponding to the background process being in state $k$. The wealth redistribution, as taking place at the ``shocks'', corresponds to the means (in self-evident notation) $w_{{\rm LL},k}$, $w_{{\rm LF},k}$, $w_{{\rm FL},k}$, and $w_{{\rm FF},k}$, again for the background process in state $k$. 
 We thus obtain
 \begin{align*}
     m'_{{\rm L},k}(t)&=\sum_{\ell=1}^dq_{\ell k}m_{{\rm L},\ell}(t)+\lambda_{{\rm L},k}\pi_k(t)+
     \gamma_k\big((w_{{\rm LL},k}-1)m_{{\rm L},k}(t)+ (I-1)w_{{\rm FL},k} m_{{\rm F},k}(t)\big), \\
     m'_{{\rm F},k}(t)&=\sum_{\ell=1}^dq_{\ell k}m_{{\rm F},\ell}(t)+\lambda_{{\rm F},k}\pi_k(t)+
     \gamma_k\big((w_{{\rm FF},k}(I-1)-1)m_{{\rm F},k}(t)+ w_{{\rm LF},k} m_{{\rm L},k}(t)\big).
 \end{align*}
 The model further simplifies if we consider the setting without modulation. In self-evident notation, we obtain ${\bs m}'(t) = \gamma(\bar A-{\mathbb I})\,{\bs m}(t)+{\bs \lambda}$, where ${\mathbb I}$ denotes a 2-dimensional identity matrix,
 \[{\bs m}(t)=\left(\begin{array}{c}
      m_{\rm L}(t)  \\
       m_{\rm F}(t)
 \end{array}\right),\:\:\:\bar A:= \left(\begin{array}{cc}
 w_{\rm LL}&(I-1)w_{\rm FL}\\
 w_{\rm LF}&(I-1)w_{\rm FF} \end{array}
 \right),\:\:\:{\bs\lambda}=\left(\begin{array}{c}
      \lambda_{\rm L}  \\
       \lambda_{\rm F}
 \end{array}\right).\]
 
 \subsubsection{Two internally homogeneous interacting subpopulations} \label{subsub} We now consider a generalization of the situation with a leader and homogeneous followers, viz.\  the situation of $I_{\rm A}$ agents of subpopulation A and $I_{\rm B}:=I-I_{\rm A}$ agents of subpopulation $B$. We let $m_{{\rm A},k}(t)$ (respectively $m_{{\rm B},k}(t)$) denote the mean wealth of an arbitrary agent from subpopulation A (respectively subpopulation B) at time $t$; the agents within each of the two subpopulations ``start symmetrically''. The arrival rates $\lambda_{{\rm A},k}$ and $\lambda_{{\rm B},k}$ are defined in the evident manner, and so are the means  $w_{{\rm AA},k}$, $w_{{\rm AB},k}$, $w_{{\rm BA},k}$, and $w_{{\rm BB},k}$. Using the same reasoning as above, we obtain
 \begin{align*}
     m'_{{\rm A},k}(t)&=\sum_{\ell=1}^dq_{\ell k}m_{{\rm A},\ell}(t)+\lambda_{{\rm A},k}\pi_k(t)+
     \gamma_k\big((I_{\rm A}w_{{\rm AA},k}-1)m_{{\rm A},k}(t)+ I_{\rm B}w_{{\rm BA},k} m_{{\rm B},k}(t)\big), \\
     m'_{{\rm B},k}(t)&=\sum_{\ell=1}^dq_{\ell k}m_{{\rm B},\ell}(t)+\lambda_{{\rm B},k}\pi_k(t)+
     \gamma_k\big((I_{\rm B}w_{{\rm BB},k}-1)m_{{\rm B},k}(t)+ I_{\rm A}w_{{\rm AB},k} m_{{\rm A},k}(t)\big).
 \end{align*}
  We get a further simplification in case that there is no modulation. In self-evident notation, we again obtain ${\bs m}'(t) = \gamma(\bar A-{\mathbb I})\,{\bs m}(t)+{\bs \lambda}$, but now with
 \[{\bs m}(t)=\left(\begin{array}{c}
      m_{\rm A}(t)  \\
       m_{\rm B}(t)
 \end{array}\right),\:\:\:\bar A:= \left(\begin{array}{cc}
 I_{\rm A}w_{\rm AA}&I_{\rm B}w_{\rm BA}\\
 I_{\rm A}w_{\rm AB}& I_{\rm B}w_{\rm BB}\end{array}
 \right),\:\:\:{\bs\lambda}=\left(\begin{array}{c}
      \lambda_{\rm A}  \\
       \lambda_{\rm B}
 \end{array}\right).\]

\section{Derivation of second moments}\label{sec:SM}
In this section we focus on characterizing the second moments pertaining to the vector ${\bs M}(t)$. The techniques relied upon resemble those used in the previous section to find the first moments. In particular, the solution again amounts to solving a system of linear differential equations.

\subsection{Differential equations for second moments}\label{subsec:2M}
Concretely, our aim is to  provide recipes to evaluate the reduced second moments of $M_i(t)$, i.e.,  
\[v_{iik}(t) := {\mathbb E}[M_i(t)(M_i(t)-1)\mathbbm{1}_{\{X(t)=k\}}]=\frac{\partial^2 f_k(\boldsymbol{z},t)}{\partial z_i^2},\]
as well as the mixed second moments of $M_i(t)$ and $M_{i'}(t)$, i.e., for $i\not =i'$,
\[v_{ii'k}(t) := {\mathbb E}[M_i(t)\,M_{i'}(t)\mathbbm{1}_{\{X(t)=k\}}]=\frac{\partial^2 f_k(\boldsymbol{z},t)}{\partial z_i\,\partial z_{i'}}.\]
Again the idea is to set up a system of coupled linear differential equations.  In these differential equations now both the transient state probabilities $\pi_k(t)$ and the transient first moments $m_{ik}(t)$ feature.
With these objects at our proposal (recalling in particular that an expression for $m_{ik}(t)$ was identified in the previous section), we can determine the corresponding variances and covariances in the evident manner. 

In the derivation, we use the identity, for $i,i'=1,\ldots,I$,
\begin{align*}\frac{\partial^2 f_k({\bs h}_{k}(\boldsymbol{z}),t)}{\partial z_i \,\partial z_{i'}}&=
\frac{\partial }{\partial z_i}\left(
\sum_{j'=1}^I \frac{\partial f_k({\bs x},t) }{\partial x_{j' }}\Big|_{{\bs x}={\bs h}_k({\bs z})} \frac{\partial (h_k({\bs z}))_{j' }}{\partial z_{i'}}\right)\\
&= \sum_{j=1}^I \sum_{j'=1}^I \frac{\partial^2 f_k({\bs x},t) }{\partial x_j\,\partial x_{j'}}\Big|_{{\bs x}={\bs h}_k({\bs z})} \frac{\partial (h_k({\bs z}))_j}{\partial z_i}\frac{\partial (h_k({\bs z}))_{j'}}{\partial z_{i'}}\\
&\:\:\:\:+\:\sum_{j' =1}^I \frac{\partial f_k({\bs x},t) }{\partial x_{j' }}\Big|_{{\bs x}={\bs h}_k({\bs z})} \frac{\partial^2 (h_k({\bs z}))_{j' }}{\partial z_i \,\partial z_{i'}},
\end{align*}
which is based on standard differentiation rules. 
We thus obtain, by differentiating \eqref{eq_PGFDifToTime} with respect to $z_i$ and $z_{i'}$ and inserting ${\bs z}={\bs 1}$, 
\begin{align*}v_{ii'k}'(t)&= \sum_{\ell=1}^d q_{\ell k} v_{ii'\ell}(t)+
\mathbbm{1}_{\{i\not=i'\}}\sum_{j:i,i'\in S_j}\lambda_{jk} \,\pi_{k}(t)+\sum_{j:i\in S_j}\lambda_{jk}m_{i'k}(t) + \sum_{j:i'\in S_j} \lambda_{jk}m_{ik}(t)\:+\\
&\:\:\:\gamma_k\left(
\sum_{j=1}^I \sum_{j'=1}^I v_{jj'k}(t)\,w_{jik}\,w_{j'i'k}+\sum_{j=1}^I m_{jk}(t)\,w^{(2)}_{jii'k}
- v_{ii'k}(t)\right);\end{align*}
here
\[w^{(2)}_{jii'k}:= \frac{\partial^2 (h_k({\bs z}))_j}{\partial z_i \,\partial z_{i'}}\Big|_{{\bs z}={\bs 1}},\]
which equals ${\mathbb E}[W_{jik}(W_{jik}-1)]$ if $i=i'$ and ${\mathbb E} [W_{jik}W_{ji'k}]$ otherwise.

\subsection{Special case: two subpopulations}
We consider the situation of Section \ref{subsub} and derive the differential equations for the (reduced) second moments in case there is no modulation. Five quantities are to be determined:
\begin{align*}
    v_{\rm AA}(t) &:= \mbox{reduced 2nd moment of arbitrary agent in population A,}\\
    v_{\rm BB}(t) &:= \mbox{reduced 2nd moment of arbitrary agent in population B,}\\
    v_{\rm AA'}(t) &:= \mbox{mixed 2nd moment of two arbitrary distinct agents in population A,}\\
    v_{\rm BB'}(t) &:= \mbox{mixed 2nd moment of two arbitrary distinct agents in population B,}\\
    v_{\rm AB}(t) &:= \mbox{mixed 2nd moment of two arbitrary agents in populations A and B,}
\end{align*}
with all quantities on the right-hand side being evaluated at time $t\geqslant 0.$ The vector ${\bs v}(t)\in {\mathbb R}_+^5$ consists of the above five entries. We can write, with ${\mathbb I}$ here denoting a 5-dimensional identity matrix, 
\[{\bs v}'(t) = \gamma(\bar A-{\mathbb I})\,{\bs v}(t) + A_m\,{\bs m}(t),\]
for a suitably chosen $(5\times 5)$-matrix $\bar A$ and a suitably chosen $(5\times 2)$-matrix $A_m$. The matrix $\bar A$ is given by, with $J_x:=I_{x}(I_x-1)$ for $x\in\{{\rm A},{\rm B}\}$, 
\[\bar A=\left(\begin{array}{ccccc}
\IA(w_{\rm AA})^2&\IB(w_{\rm BA})^2&\JA(w_{\rm AA})^2&\JB(w_{\rm BA})^2&2\IA\IB w_{\rm AA}w_{\rm BA}\\
\IA(w_{\rm AB})^2&\IB(w_{\rm BB})^2&\JA(w_{\rm AB})^2&\JB(w_{\rm BB})^2&2\IA\IB w_{\rm AB}w_{\rm BB}\\
\IA(w_{\rm AA})^2&\IB(w_{\rm BA})^2&\JA(w_{\rm AA})^2&\JB(w_{\rm BA})^2&2\IA\IB w_{\rm AA}w_{\rm BA}\\
\IA(w_{\rm AB})^2&\IB(w_{\rm BB})^2&\JA(w_{\rm AB})^2&\JB(w_{\rm BB})^2&2\IA\IB w_{\rm AB}w_{\rm BB}\\
\IA w_{\rm AA} w_{\rm AB}&\IB w_{\rm BA} w_{\rm BB}&\JA w_{\rm AA} w_{\rm AB}&\JB  w_{\rm BA} w_{\rm BB}&\IA\IB (w_{\rm AA}w_{\rm BB} + w_{\rm AB}w_{\rm BA})\end{array}
\right)\]
and
\[A_m=
\left(\begin{array}{cc}
2\lambda_{\rm A} +\gamma \IA w_{\rm AAA}^{(2)}&\gamma \IB w_{\rm BAA}^{(2)}\\
\gamma \IA w_{\rm ABB}^{(2)}&
2\lambda_{\rm B} +\gamma \IB w_{\rm BBB}^{(2)}\\
2\lambda_{\rm A} +\gamma \IA w_{\rm AAA'}^{(2)}&\gamma \IB w_{\rm BAA'}^{(2)}\\
\gamma \IA w_{\rm ABB'}^{(2)}&
2\lambda_{\rm B} +\gamma \IB w_{\rm BBB'}^{(2)}\\
\lambda_{\rm A}+ \gamma \IA w_{\rm AAB}^{(2)}&
\lambda_{\rm B} +\gamma \IB w_{\rm BAB}^{(2)}\\
\end{array}\right).\]

\section{Application 1: wealth redistribution}\label{Sec:WD}
The model that we consider in this paper can be interpreted as a simple formalism describing an economy, providing insight into the stochastic evolution of the wealth of each of the individual agents. Indeed, the set $I$ could correspond to the agents of the economic system under study, which is fed by external inflow and in which at random times wealth redistribution occurs.  In this section we present an example of such an economic system, with a population consisting of one agent (the ``leader'') obtaining income from outside the system, and the other $I-1$ agents obtaining their income from the leader (the ``followers''). The Markovian background process records the state of the economy, in that it alternates between periods of economic growth and periods of recession. Our model allows us to quantify the distribution of the fraction of followers whose income drops below a critical threshold. In this way, we can get insight into the phenomenon of ``poverty trap'', i.e., persistent poverty for the followers (to be interpreted as the segment of the population that does not own resources, and whose income strongly depends on payments from the leader). 
The model considered is formally described as follows.

\vb

\noindent $\circ$
The background process $X(t)$ has two states, i.e., economic growth (corresponding to state 1) and recession (corresponding to state 2).

\vb

\noindent $\circ$
Regarding the income rates $\lambda_j$, we consider the situation that $S_j=\{j\}$, for $j=1,\ldots,I$. Only agent 1 (the leader) has external income: we let $\lambda_{j1}=\lambda_1$ and $\lambda_{j2}=\lambda_2$ for rates $\lambda_1$ and $\lambda_2$ such that $0<\lambda_2<\lambda_1$ (i.e., the leader has a higher income during periods of economic growth). The other agents (i.e., the followers) do not have any external income: $\lambda_{jk}=0$ for $j=2,\ldots,I.$

\vb

\noindent $\circ$ Wealth distribution takes place at a Poisson rate $\gamma_k$, with $k$ the state of the background process. The vectors $(W_{11k},\ldots,W_{1Ik})$, for $k=1,2$, are multinomially distributed with parameters $1$ and 
$(p_k, r_k,\ldots,r_k)$. Here $r_k\leqslant (1-p_k)/(I-1)$, meaning that with probability $1-p_k-(I-1)r_k\in[0,1]$ the wealth unit of agent~1 leaves the economy. 
Typically one expects $p_2>p_1$, as in periods of recession the leader will be inclined to save a larger fraction of their wealth.  We let the vectors $(W_{j1k},\ldots,W_{jIK})$, for $k=1,2$ and $j=2,\ldots,I$, be multinomially distributed with parameters $1$ and 
$(0, s_k,\ldots,s_k)$. Here $s_k\leqslant 1/(I-1)$, meaning that with probability $1-(I-1)s_k\in[0,1]$ the wealth unit of client $j$ leaves the economy. 

\vb

In this setup intentionally various symmetries are assumed, so as to keep the model as low-dimensional as possible. Evidently, more general variants can be analyzed as well. For example, we could consider instances such that, for $j=2,\ldots,I$ and $j'\not= j$, the distribution of $W_{jjk}$ differs from the distribution of $W_{jj'k}$.

\subsection{Means and variances}
\begin{figure}[ht!]
\centering
\begin{subfigure}[t]{0.45\textwidth}
\includegraphics[width=0.95\linewidth]{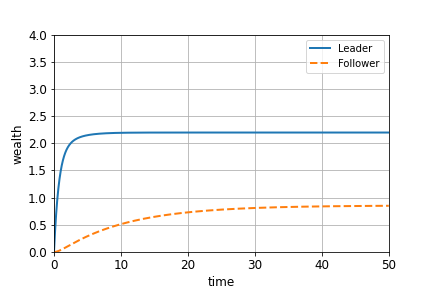} 
\caption{Transient means as numerical solution of differential equations}
\end{subfigure}\:\:\:\:
\begin{subfigure}[t]{0.45\textwidth}
\includegraphics[width=0.95\linewidth]{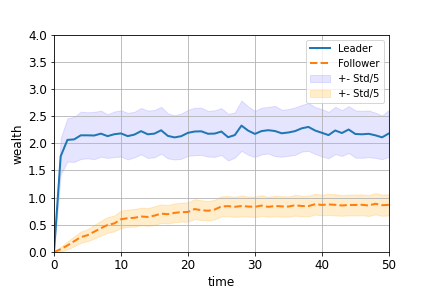} 
\caption{Transient means approximated by simulations}
\end{subfigure}\:\:\:\:

\caption{Transient means of the leader's and follower's wealth. The numerical solutions have been evaluated by applying Euler's method. For the simulated approximation we have run $2000$ simulations for a population of $I=30$. The chosen parameters are: $q_{12}=1/100$, $q_{21}=5/100$, $\lambda_1 = 3$, $\lambda_2=1$, $\gamma_1 = 2$, $\gamma_2 = 1$, $p_1 = 0.3$, $p_2 = 0.6$, $1-p_1-(I-1)r_1 = 5/100$, $1-p_2-(I-1)r_2 = 1/100$, $1-(I-1)s_1 = 5/100$ and $1-(I-1)s_2 = 1/10$.}
\label{fg:transm}
\end{figure}

It is not hard to verify that
\[g_{1k}({\bs z}) = z_1^{p_k}\left(\prod_{j=2}^I z_j\right)^{r_k},\:\:\:\:
g_{jk}({\bs z})= \left(\prod_{j=2}^I z_j\right)^{s_k},\]
thus also defining $h_k({\bs z})$. For ease, assume that the background process is in stationarity at time $0$.
We thus obtain the following differential equations for the transient means (in self-evident notation):
\begin{align*}
    m_{{\rm L},k}'(t)&= \sum_{\ell=1}^2q_{\ell k}m_{{\rm L},\ell}(t)+\lambda_{k}\pi_k +\gamma_k(p_k-1)\,m_{{\rm L},k}(t),\\
    m_{{\rm F},k}'(t)&= \sum_{\ell=1}^2q_{\ell k}m_{{\rm F},\ell}(t)+\gamma_k\big(((I-1)s_k-1)m_{{\rm F},k}(t)+r_k m_{{\rm L},k}(t)\big),
\end{align*}
using the results found in Section \ref{ssec:LHF}.
{The transient means as described by the differential equations above indeed nicely match simulation results as shown in Figure \ref{fg:transm}.}
Focusing on stationarity, we obtain that
\[\left(\begin{array}{c}m_{{\rm L},1}\\m_{{\rm L},2}\end{array}\right) = -  \left(\begin{array}{cc}-q_1+\gamma_1(p_1-1)&q_2\\q_1&-q_2+\gamma_2(p_2-1)\end{array}\right)^{-1}\left(\begin{array}{c}\lambda_1\pi_1\\ \lambda_2\pi_2\end{array}\right)\]
and
\begin{align*}
    \left(\begin{array}{c}m_{{\rm F},1}\\m_{{\rm F},2}\end{array}\right) &= -   \left(\begin{array}{cc}-q_1+\gamma_1((I-1)s_1-1)&q_2\\q_1&-q_2+\gamma_2((I-1)s_2-1)\end{array}\right)^{-1}\cdot\\
    &\:\:\:\hspace{3cm}\left(\begin{array}{cc}\gamma_1r_1&0\\0&\gamma_2r_2\end{array}\right)\left(\begin{array}{c}m_{{\rm L},1}\\m_{{\rm L},2}\end{array}\right).
\end{align*}
Using Proposition \ref{pr:stab} we can determine under what condition these stationary means are well defined. 
A similar system of equations can be set up to determine the corresponding stationary (reduced) second moments, in self-evident notation denoted by
\[(v_{{\rm LL},1},v_{{\rm LL},2}, v_{{\rm FF},1} , v_{{\rm FF},2}, v_{{\rm FF'},1}, v_{{\rm FF'},2}, 
v_{{\rm LF},1}, v_{{\rm LF},2});\]
we do not present the expressions here.

\subsection{Poverty trap for single follower}
\begin{figure}[htbp]
\centering
\includegraphics[width=0.45\linewidth]{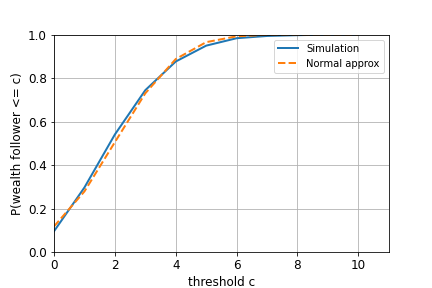} 
\caption{\label{fg:thresh} Cumulative distribution function of the follower's wealth, approximated through simulations and the normal approximation. For the simulated approximation we have run $5000$ simulations for a population of $I=30$. The chosen parameters are: $q_{12}=1/100$, $q_{21}=5/100$, $\lambda_1 = 10$, $\lambda_2=6$, $\gamma_1 = 4$, $\gamma_2 = 2$, $p_1 = 0.2$, $p_2 = 0.4$, $1-p_1-(I-1)r_1 = 5/100$, $1-p_2-(I-1)r_2 = 1/100$, $1-(I-1)s_1 = 3/100$ and $1-(I-1)s_2 = 7/100$.}
\end{figure}

Let $f$ be the probability of the stationary wealth of an arbitrary follower, denoted by $M_{j}$ for some follower $j=2,\ldots,I$, being below some critical threshold $c$; this $f$ can, in our stylized context, be considered as the probability of a follower ending up in the poverty trap. 
Using a straightforward normal approximation, we conclude
\[f = {\mathbb P}(M_{j}\leqslant c) \approx f_N:=\Phi\left(\frac{c + 1/2 -m_{\rm F}}{\sqrt{v^\circ_{\rm FF}}}\right),\]
where 
\[m_{\rm F}:=  m_{{\rm F},1}+ m_{{\rm F},2},\:\:\:\:\:v_{\rm FF}:=  v_{{\rm FF},1}+ v_{{\rm FF},2}\]
and
 $v^\circ_{{\rm FF}}:=v_{{\rm FF}}+m_{{\rm F}}-m_{{\rm F}}^2$;
as usual, $\Phi(\cdot)$ denotes the cumulative distribution function of a standard normal random variable. To verify the accuracy of our findings, we have also estimated ${\mathbb P}(M_{j}\leqslant c)$ by performing 5000 independent simulation runs. The resulting distribution is indeed highly similar to the normal approximation, as shown in Figure \ref{fg:thresh}.

Note that the followers do not operate independently, as they react to a common background process. As a consequence, the number of followers ending up in the state of a poverty trap, say $B$, is {\it not} binomially distributed (even though all of them experience the same probability $f$ of doing so). In the remainder of this section we point out how to derive a good proxy for the distribution of $B$. 

\subsection{Poverty trap for full follower population}
\begin{figure}[htbp]
\centering
\includegraphics[width=0.45\linewidth]{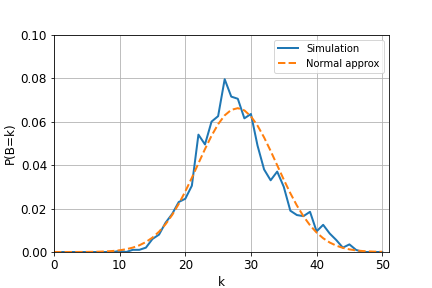} 
\caption{\label{fg:totresh} Density function of the amount of follower's with wealth less or equal to 1, evaluated through simulations and the normal approximation. For the simulated approximation we have run $2000$ simulations for a population of $I=50$. The chosen parameters are: $q_{12}=1/100$, $q_{21}=5/100$, $\lambda_1 = 10$, $\lambda_2=6$, $\gamma_1 = 4$, $\gamma_2 = 2$, $p_1 = 0.2$, $p_2 = 0.4$, $1-p_1-(I-1)r_1 = 5/100$, $1-p_2-(I-1)r_2 = 1/100$, $1-(I-1)s_1 = 3/100$ and $1-(I-1)s_2 = 7/100$.}
\end{figure}

With $B_j$ the indicator function of $M_j\leqslant c$, for $j=2,\ldots,I$, we are interested in approximating the distribution of
\[B:=\sum_{j=2}^I B_j.\]
We propose a normal approximation, which relies on the central limit theorem, entailing that it will be particularly accurate as $I$ grows. Define, for $j,j'=2,\ldots,I$ such that $j\not = j'$,
\[f' := {\mathbb P}(M_{j}\leqslant c, M_{j'}\leqslant c).\]
This probability can be approximated by its Gaussian counterpart. Applying continuity correction, 
\[f'\approx f'_N:={\mathbb P}\left(M^\circ_{j}\leqslant c +\frac{1}{2}, M^\circ_{j'}\leqslant c +\frac{1}{2}\right),\]
with $(M^\circ_{j}, M^\circ_{j'})$ bivariate normal with mean $(m_{\rm F},m_{\rm F})$ and covariance matrix
\[\Sigma = \left(\begin{array}{cc}v^\circ_{{\rm FF}}&v^\circ_{{\rm FF'}}\\v^\circ_{{\rm FF'}}&v^\circ_{{\rm FF}}\end{array}\right),\]
where  $v^\circ_{{\rm FF}}:=v_{{\rm FF}}+m_{{\rm F}}-m_{{\rm F}}^2$ (as before) and $v^\circ_{{\rm FF'}}:=v_{{\rm FF'}}-m_{{\rm F}}^2$.    
Notice that there are powerful numerical techniques to accurately evaluate bivariate normal probabilities; see for instance~\cite{COX}.
It now follows that
\begin{align*}{\mathbb E}\,B&\approx\mu_B:=(I-1)f_N,\\{\mathbb V}{\rm ar}\,B&\approx \sigma^2_B:=(I-1)f_N(1-f_N)+(I-1)(I-2)(f'_N-f_N^2).\end{align*}
Applying a standard continuity correction, this gives rise to the following approximation: for $k=1,\ldots,I-1$,
\[{\mathbb P}(B=k) \approx \Phi\left(\frac{k+\frac{1}{2}-\mu_B}{\sigma_B}\right)- \Phi\left(\frac{k-\frac{1}{2}-\mu_B}{\sigma_B}\right).\]
As shown in Figure \ref{fg:totresh}, the distribution following from the normal approximation matches the simulation-based approximation quite well.

\section{Application 2: opinion dynamics}\label{sec:OD}
In the field of opinion dynamics one studies, predominantly based on mathematical  models, the evolution of opinions in a collection of agents. Arguably the most basic, yet meaningful model was proposed in \cite{DEG}, considering a set of $n$ agents, each of them having an opinion on a particular subject. At every discrete time instant each agent updates her opinion based on the other agents' opinions. More specifically, the opinion of an agent's opinion at a certain point in time is a weighted sum of all agents' opinion at the previous time instant, with the weights summing to $1$.  To make the model more realistic, various extensions have been developed. In this context, one of the important contributions is \cite{FRI}, in which agents are assumed to have an internal opinion as well as an {expressed} opinion. In modelling terms, it means that the framework of \cite{FRI} extends the one of \cite{DEG} by dropping the assumption that the weights sum to $1$ and by allowing (at any point in time) an external contribution to the opinion vector. 

The goal of this section is to show that one can use  the theory developed in the present paper to make opinion dynamics models considerably more realistic.
Compared to existing frameworks on opinion dynamics, two significant improvements can be made. In the first place, most of the existing models on opinion dynamics that allows for analytical investigation are of a deterministic nature, and as such do not incorporate the random effects present in a population of agents influencing each other's opinions; these can be included relying on our modelling framework. In the second place, it allows us, through the modulation mechanism, to incorporate randomly evolving external effects that have impact on the agents' opinions. {Examples of external effects could relate to the current level of economic prosperity, or to the degree of access to digital communication, or to the geographic area agents live in; these variables change over time (where the natural timescale could correspond to years) which result in different opinion dynamics.}

Let us consider the first advantage, i.e., the option of introducing stochasticity, in more detail. The framework that we developed can in fact be interpreted as a stochastic generalization of the model studied in \cite{FRI}, which (in its most elementary form) can be summarized as follows.
Let $Y_t$ a vector of opinions, $W_t$ a matrix that describes the effects of each opinion held at time $t-1$ on the opinions held at time $t$, $X_t$ a matrix of scores on exogenous variables, $B_t$ a vector of coefficients giving the effects of each of the exogenous variables, $\alpha_t$ a scalar weight corresponding to the endogenous conditions and $\beta_t$  a scalar weight corresponding to the exogenous conditions. 
Then the recurrence relation that defines the model in \cite{FRI} is given, for $t\in {\mathbb N}$, by
\[Y_t = \alpha_t W_t Y_{t-1} + \beta_t X_t B_t.\]
The relation with our model can be seen immediately by comparing this recurrence relation with the differential equation of the transient means in Proposition \ref{PRO:MDE}. {Indeed, just as in \cite{FRI}, the opinion of an agent is construed as a linear combination of its own and the others' opinions before the time of redistribution.} An important additional advantage of our approach is the option to explicitly characterize all (reduced) moments of the opinion vector through systems of differential equations, by repeated differentiation of the relation featuring in Proposition \ref{PRO:FDE}. To our best knowledge, this is a novel analytical tool that was not provided in other stochastic generalizations of the Friedkin-Johnsen model \cite{FRI}.

\subsection{Opinion dynamics with Markovian background process}
In this example we consider opinion dynamics in a two-group population (groups A and B) where individuals of one group -- say group A --  have an opinion update mechanism that alternates between a ``normal mode'' and an ``adapted mode''. This framework may be interpreted as a stylized model for groups of individuals occasionally visiting larger events as conferences, demonstrations, or political gatherings.
Formally we describe the model as follows.

\vb

\noindent $\circ$ The background process $X(t)$ has two states, i.e., one corresponding to the normal mode  (corresponding to state 1) and one to the adapted mode (corresponding to state 2).

   \vb

\noindent $\circ$ We consider the situation that individuals do not automatically {increase the value of} their opinions when the background process is in the normal mode. Since group~B agents are not affected by the background process, we thus have $\lambda_{\rm A,1}=\lambda_{\rm B,1}=\lambda_{\rm B,2}=0$. Additionally, we let group A agents automatically strengthen their opinion with rate $\lambda_{\rm A,2} \geqslant 0$ in the adapted mode, i.e., they only strengthen their opinion as a result of external influence by attending conference, political meetings, etc.

\vb

\noindent $\circ$ We let the Poisson rate of opinion redistribution $\gamma > 0$ be unaffected by the background process. Since group A separates itself from group B when the background process is in the adapted state, we assume that group A and B individuals do not affect each other's opinion in that state. In addition, in this example we allow group A agents to be ``overenthusiastic'' in the adapted state, resulting in distributing more opinion units than they originally have.
    
    We first consider the situation that the background process is in state $1$.
    Let there be $I_{\rm A}$ agents in group A, and $I_{\rm B}$ in group B. For $j\in\{1,\ldots, I_{\rm A}\}$ (i.e., corresponding to agents in group A), we let the vector
    \[(W_{j,1,1}, \ldots, W_{j, I_{\rm A},1},W_{j, I_{\rm A} + 1, 1}, \ldots, W_{j ,I_{\rm A}+ I_{\rm B},1})\] be multinomially distributed with parameters $1$ and $(p_{\rm AA,1},\ldots,p_{\rm AA,1},p_{\rm AB,1},\ldots,p_{\rm AB,1})$; for $j\in\{I_{\rm A} +1 ,\ldots, I_{\rm A} + I_{\rm B}\}$ (i.e., corresponding to agents in group B) these parameters are $1$ and $(p_{\rm BA,1},\ldots,p_{\rm BA,1},p_{\rm BB,1},\ldots,p_{\rm BB,1})$. Here we consider the situation that \[I_{\rm A} \: p_{\rm AA,1} + I_{\rm B} \: p_{\rm AB,1} = 1\:\:\:\mbox{ and }\:\:\: I_{\rm A} \: p_{\rm BA,1} + I_{\rm B} \: p_{\rm BB,1} = 1,\] effectively meaning that ``opinion mass cannot leak away from the system'' when the background process is in state 1.
    
    When the background process is in state $2$, we let for $j\in\{I_{\rm A} +1 ,\ldots, I_{\rm A} + I_{\rm B}\}$ (i.e., for agents of group B) \[(W_{j,1,2}, \ldots, W_{j, I_{\rm A},2},W_{j, I_{\rm A} + 1, 2}, \ldots, W_{j ,I_{\rm A}+ I_{\rm B},2})\]  be multinomially distributed with parameters $1$ and $(0,\ldots,0,p_{\rm BB,2},\ldots,p_{\rm BB,2})$. Here, group B opinion does not leak away: $I_B \: p_{\rm BB,2} = 1$. For $j\in\{1,\ldots, I_{\rm A}\}$ (i.e., for agents of group A) we introduce an additional mechanism: with probability $0 \leqslant \alpha \leqslant 1$ the unit of opinion to be distributed becomes two units and with probability $1-\alpha$ it remains one unit. Afterwards, each unit is multinomially distributed with parameters $1$ and $(p_{\rm AA,2},\ldots,p_{\rm AA,2},0,\ldots,0)$, with $I_{\rm A} \: p_{\rm AA,2} = 1$.

\vb

We would like to stress that the mechanism described above is just an example of one specific type of opinion dynamics we can handle. Various alternative models can be dealt similarly (more than two groups, more than two background states, the ``$\alpha$-jumps'' corresponding with multiplying the ``opinion unit'' with a number different from two, etc.). 

The model proposed allows groups of individuals of group A to ``create opinion'' when the background state is 2 (by $\lambda_{\rm A,2} > 0$), and to multiply each of their existing opinion units by $2$ with probability $\alpha$. Since ``opinion mass'' does not leave the population in this setup, we typically expect opinions to grow unboundedly in time when $\lambda_{\rm A,2} > 0$ or $\alpha > 0$. In the next section we show that the aforementioned is indeed true for the means of the agent's opinion.

\subsection{Means and variances}
For ease, we assume that the background process is in stationarity at time 0. With the results found in Section \ref{subsub} and the notation as in Proposition~\ref{PRO:MDE}, we have the following collection of differential equations describing the transient means: \[{\bs m}'(t) = A\,{\bs m}(t)+\Lambda\,{\bs \pi},\] with ${\bs m}(t) = (m_{\rm A, 1}(t), m_{\rm A, 2}(t), m_{\rm B, 1}(t), m_{\rm B, 2}(t))^{\top}$,

\[A=\left(\begin{array}{cccc}
-q_{12} + \gamma (\IA \, p_{\rm {AA},1} - 1) & q_{21} & \gamma \, \IB \, p_{\rm {BA},1} & 0\\
q_{12} & -q_{21} + \gamma \, \alpha & 0 & 0\\
\gamma \, \IA \, p_{\rm {AB},1} & 0 & -q_{12} + \gamma \, (\IB \, p_{\rm {BB},1} -1)&q_{21}\\
0 & 0 & q_{12} & -q_{21}
\end{array}
\right)\]
and $\Lambda = {\rm diag}\{0, \lambda_{\rm A,2},0,0\}$. For the transient second moments, similar differential equations can be derived as well, relying on the results presented in Section \ref{subsec:2M}.

It is straightforward to verify that a necessary condition for the existence of a steady state is $\alpha =0$ and $\lambda_{\rm A,2}=0$, as anticipated in the previous section. {Under these conditions, opinions are only redistributed since no additional opinion is ``created'' and, by assumption, no opinion is ``lost''.} The steady state of this model is, by solving $A \, \boldsymbol{m} = \boldsymbol{0}$, proportional to the vector 
\begin{equation}\label{eq:sse}\left(1,\:\: \frac{q_{12}}{q_{21}}, \:\: \frac{p_{\rm {AB},1}}{p_{\rm {BA},1}},\:\: \frac{p_{\rm {AB},1} \, q_{12}}{p_{\rm {BA},1} \, q_{21}}\right)^{\top}.\end{equation}
This vector should be normalized such that the total ``opinion mass'' equals its initial value (i.e., ${\bs 1}^\top {\bs M}(0)$). 
Indeed, if the setup of the system corresponds to spending  longer periods of time in state $2$, or equivalently $q_{12} > q_{21}$, then we expect $m_{\rm A,2} > m_{\rm A,1}$. Also, when group $\rm A$ agents distribute their opinions with a higher probability to group $\rm B$ agents than the probability of group $\rm B$ agents distributing their opinion to group $\rm A$ agents, thus $p_{\rm {AB},1} > p_{\rm {BA},1}$, we expect $m_{{\rm B},1} > m_{{\rm A},1}$. Both properties are in line with the steady-state expression \eqref{eq:sse}. 

When $\alpha=0$ and $\lambda_{\rm A,2}=0$, the adapted state can be interpreted as a temporary interruption of the opinion formation process between groups $\rm A$ and $\rm B$; in the adapted state, agents of both groups do not distribute their opinion to other group's agents. This interruption should not affect the steady state of the system. Indeed, the steady states of the groups are related through $m_{\rm A} = ({p_{\rm {BA},1}}/{p_{\rm {AB},1}}) \, m_{\rm B}$ which does not involve the background process parameters $q_{12}$ and $q_{21}$.

\begin{figure}[ht!]
\centering
\begin{subfigure}[t]{0.45\textwidth}
\includegraphics[width=0.95\linewidth]{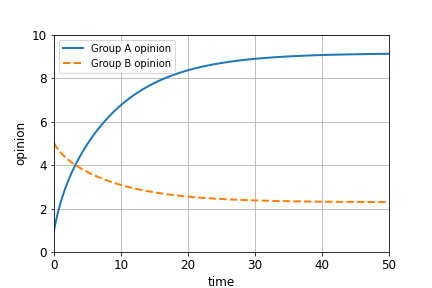} 
\caption{Transient means as numerical solution of differential equations}
\end{subfigure}\:\:\:\:
\begin{subfigure}[t]{0.45\textwidth}
\includegraphics[width=0.95\linewidth]{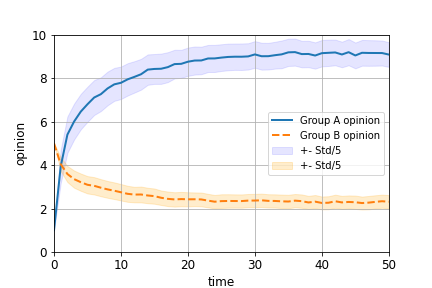} 
\caption{Transient means approximated by simulations}
\end{subfigure}\:\:\:\:

\caption{Transient means of group {\rm A}'s and {\rm B}'s opinion. The numerical solutions have been evaluated by applying Euler's method. For the simulated approximation we have run $2000$ simulations. The group sizes are $I_{\rm A} = 10$ and $I_{\rm B} = 30$. The initial opinions of group {\rm A} and {\rm B} agents are $1$ and~$5$ respectively. {The aggregate opinion of all agents, $I_{\rm A} m_{\rm A}(t) + I_{\rm B} m_{\rm B}(t)$, remains constant and  equals $I_{\rm A} {m}_{\rm A}(0) + I_{\rm B} {m}_{\rm B}(0) = 160$, as we have enforced ``conservation of opinion''.} The chosen parameters are: $q_{12}=3/10$, $q_{21}=2/10$, $\lambda_{\rm A, 2} = 0$, $\alpha = 0$, $\gamma = 5/8$, $I_{\rm A}\: p_{\rm {AA},1} / I_{\rm B}\: p_{\rm {AB},1} = 2$ and $I_{\rm A}\: p_{\rm {BA},1} / I_{\rm B}\: p_{\rm {BB},1} = 4/5$. The abscissa $\omega$ is $0$.}
\label{fg:opinion_numsim}
\end{figure}

Figure \ref{fg:opinion_numsim} shows that the simulation-based approximation of the means aligns with the  numerical solution of the differential equation of the means. In this setup the abscissa $\omega$ of the matrix $A$ is equal to $0$, thus not satisfying the stability condition of Proposition \ref{pr:stab}. Still, the transient means converge to a steady state; bear in mind that the condition of Proposition \ref{pr:stab} is a {\it sufficient} condition, hence there can be stability even though $\omega<0$ is not fulfilled. In Section \ref{ss:unbounded} we provide an example where there is no stability while the abscissa $\omega$ is~$0$. {In the setup of Figure \ref{fg:opinion_numsim}, the agents of group $\rm A$ tend to give more attention to the topic of interest than the agents of group $\rm B$, in the sense that $I_{\rm A} \,p_{\rm {AA},1} + I_{\rm B}\, p_{\rm {BA},1} = 5/3$  while $I_{\rm A} \,p_{\rm {AB},1} + I_{\rm B} \,p_{\rm {BB},1} = 8/9$. This ``unequal attentiveness'' typically results in polarization, as discussed in \cite{CHA}.}

\begin{figure}[htbp]
\centering
\includegraphics[width=0.45\linewidth]{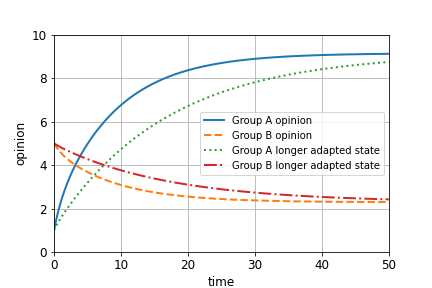} 
\caption{\label{fg:longer_adapt} Effect of a longer period of time in the adapted state on the transient means. The transient means have been evaluated by applying Euler's method. The longer period of time in the adapted state is obtained by increasing $q_{12}$ (i.e., we now set $q_{12}=1$, versus $q_{12}=3/10$ in the base case). All the other parameters are identical to the parameters used in Figure~\ref{fg:opinion_numsim}. The abscissa $\omega$ is $0$ in both situations.}
\end{figure}

In Figure \ref{fg:longer_adapt} we show the effect of different background parameters $q_{12}$ and $q_{21}$. The curve of the transient means corresponding to spending a longer fraction of time in the adapted state (i.e., the one with the higher $q_{12}$) moves more slowly to the steady state.

\subsection{Unbounded opinion mass and relative opinions} \label{ss:unbounded}
\begin{figure}[ht!]
\centering
\begin{subfigure}[t]{0.45\textwidth}
\includegraphics[width=0.95\linewidth]{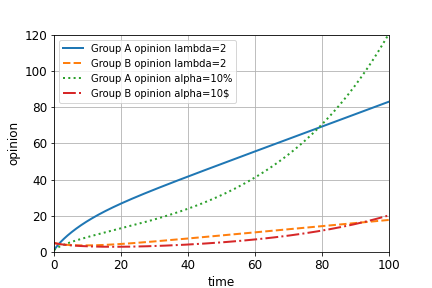} 
\caption{Transient means ``growing'' without limits}
\label{fg:opinion_absrelA}
\end{subfigure}\:\:\:\:
\begin{subfigure}[t]{0.45\textwidth}
\includegraphics[width=0.95\linewidth]{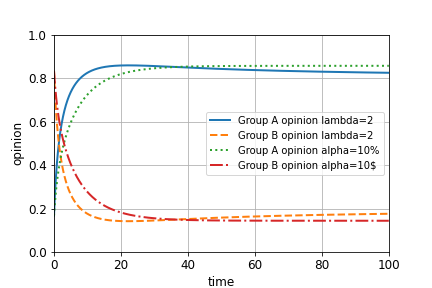} 
\caption{Transient means interpreted as relative opinions}
\label{fg:opinion_absrelB}
\end{subfigure}\:\:\:\:

\caption{Unbounded transient means become bounded when interpreted as relative opinions. The numerical solutions have been evaluated by applying Euler's method. All the parameters chosen are identical to the setup as in Figure \ref{fg:opinion_numsim} except for $\lambda_{\rm A, 2}$ and $\alpha$. The abscissa $\omega$ is $0$ in the setup with $\lambda_{\rm A, 2} = 2$. The abscissa $\omega$ is $0.027$ in the setup with $\alpha = 10\%$ .}
\label{fg:opinion_absrel}
\end{figure}

We now consider the setup where $\alpha$ or $\lambda_{{\rm A},2}$ is assumed to be strictly positive. As reasoned earlier, in this situation we expect the total ``opinion mass'' to grow beyond any bound. In Figure \ref{fg:opinion_absrelA} we show the evolution of the transient means in the situations $(\alpha, \lambda_{{\rm A},2}) = (0, 2)$ and $(\alpha, \lambda_{{\rm A},2}) = (1/10, 0)$. 

It is noted that, in principle, unbounded opinions (in our case due to $\alpha > 0$ or $\lambda_{{\rm A},2} > 0$) are an outcome of the model that cannot supported by empirical data. 

However, various techniques have been proposed in order to avert the phenomenon of unbounded opinions while maintaining the mentioned interaction mechanism, as explained in  \cite{FLA17}. One of the ideas is to interpret the modeled opinions as \textit{relative} {opinions} \cite{CHA}. This means that the absolute values of the entries of the opinion vector have no meaning, but that interpretation is only given to the fractions between these entries.  Concretely, the $i$-th entry divided by the $j$-th entry reflects the relative opinion of agent $i$ with respect to agent $j$.

In Figure \ref{fg:opinion_absrelB} it is shown how the unbounded growth of opinion mass in Figure \ref{fg:opinion_absrelA} changes into a stabilizing growth when considering the relative perspective. This illustrates that, when dealing with instances of the stochastic model where opinion mass may grow unboundedly, one may consider working with relative opinions instead. With this interpretation, one can still obtain insight into the effects that the specific interaction mechanism between the agents has  on the opinion formation process of the entire population. Further empirical and technical arguments in favor of working with relative opinions, rather than their absolute counterparts, are discussed in great detail in \cite{CHA}.

\section{Application 3: file storage systems}\label{Sec:SS}
Consider a system in which users generate files. For safety reasons, the files that have been saved at the clients' locations are periodically copied to a central storage location where a backup is made. The frequency of copying the clients' files has to be sufficiently high to make sure that a relatively small amount has no centrally stored duplicate. In this section we discuss a sequence of models describing the dynamics of such file storage systems, staring with the most rudimentary variant. 
For more background on various aspects of data storage networks, we refer to e.g.\ 
\cite{PES}. 

\subsection{Basic variant}
A stylized model that describes the dynamics of this system is the following. 
Let $M_1(t)$ record the number of files at the clients' devices that have not been copied to the central storage unit by time $t\geqslant 0$, and let $M_2(t)$ the number of files at the storage unit at time $t\geqslant 0$. Let $\lambda$ be the Poisson rate at which the aggregate client population generates files. Let the times between subsequent backups be exponentially distributed with mean $\gamma^{-1}.$ In this system without modulation, we observe that $W_{12}\equiv W_{22}\equiv 1$, and $W_{11}\equiv W_{21}\equiv 0$.

It can be verified from the results of Section \ref{sec:FM} that
\[m'_1(t) = \lambda - \gamma m_1(t),\:\:\:\:m_2'(t) = \gamma m_1(t).\]
It is not hard to verify that  the first moments can be computed explicitly: if the system starts empty, for $t\geqslant 0$,
\[m_1(t) = \frac{\lambda}{\gamma}\big(1-e^{-\gamma t}\big),\:\:\:\:m_2(t) = {\lambda t}  - \frac{\lambda}{\gamma}\big(1-e^{-\gamma t}\big).\]
These expressions can be used to determine the optimal backup rate $\gamma^\star$, for instance by considering a cost function that encompasses the per-backup cost (with proportionality constant $\kappa_{\rm B}>0$) and the  cost associated with the number of files that have not been copied yet (with proportionality constant $\kappa_{\rm NC}>0$). Concretely, for given time horizon $t\geqslant 0$, this leads to the optimization problem
\[\min_{\gamma>0}F_1(\gamma)+F_2(\gamma),\:\:\:\:\:\mbox{with}\:\:F_1(\gamma):= \gamma t\, \kappa_{\rm B} ,\:\:\:F_2(\gamma):= \frac{\lambda}{\gamma}\big(1-e^{-\gamma t}\big) 
\kappa_{\rm NC},\]
with the objective function $F(\gamma):=F_1(\gamma)+F_2(\gamma)$ being convex.
It takes some calculus to verify that if $F'(0)\geqslant 0$, or (equivalently) $2\kappa_{\rm B}\geqslant \lambda t \,\kappa_{\rm NC}$, then $\gamma=\gamma^\star=0$ is optimal, which in practical terms means that one should very frequently update; if on the other hand $F'(0)< 0$, or (equivalently) $2\kappa_{\rm B}<\lambda t \,\kappa_{\rm NC}$, there is a strictly positive optimal update rate $\gamma^\star$. In the latter case, $\gamma^\star$ cannot be computed in closed form, but that can be easily numerically evaluated; see also Figure \ref{F1}

{\footnotesize
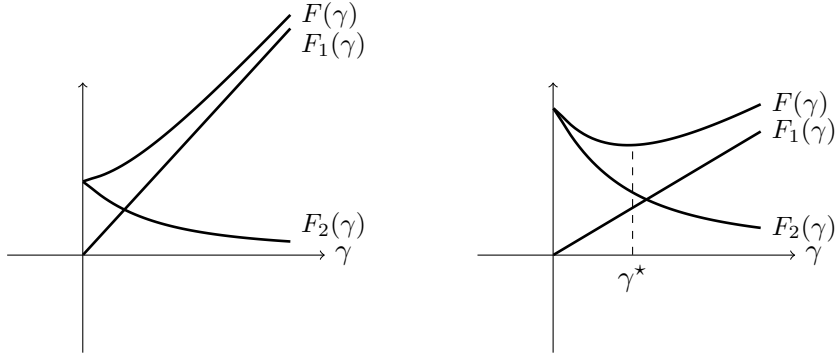
\begin{figure}
\hspace{1.2cm}
    \begin{tikzpicture}[scale=0.9950]
%\begin{axis}[xmax=9,ymax=9, samples=50]
      \draw[->] (-1,0) -- (3.2,0) node[right] {$\gamma$};
      \draw[->] (0,-1.3) -- (0,2.3);
   %   \draw (1,0) -- (1,-0.2) node[below] {$\psi(\beta)$}; 
   %   \draw (1,0) -- (1, 0.95);
%      \draw (1,0.95) -- (-0.2,0.95) node[left] {$\beta$}; 
      \draw[scale=0.5,domain=0.01:5.49,smooth,variable=\x,black,line width=1pt] plot ({\x},{1.1*\x}); %node[right] {\small $F_1(\gamma)$}; 
      \draw (2.75,2.8) node[right] {\small $F_1(\gamma)$};
       \draw (2.75,0.4) node[right] {\small $F_2(\gamma)$};
       \draw[scale=0.5,domain=0.01:5.49,smooth,variable=\x,black,line width=1pt] plot ({\x},{2*(1-exp(-\x))/\x});
       %node[right] {\small $F_2(\gamma)$};
       \draw[scale=0.5,domain=0.01:5.49,smooth,variable=\x,black,line width=1pt] plot ({\x},{2*(1-exp(-\x))/\x + 1.1*\x}) node[right] {\small $F(\gamma)$};
      %     \draw[scale=0.5,domain=-1.92:1,smooth,variable=\x,blue] plot ({\x},{1/3*1/(2+\x)});
  %    \end{axis}
    \end{tikzpicture}\hspace{1.2cm}
    \begin{tikzpicture}[scale=0.9950]
%\begin{axis}[xmax=9,ymax=9, samples=50]
      \draw[->] (-1,0) -- (3.2,0) node[right] {$\gamma$};
      \draw[->] (0,-1.3) -- (0,2.3);
   %   \draw (1,0) -- (1,-0.2) node[below] {$\psi(\beta)$}; 
   %   \draw (1,0) -- (1, 0.95);
%      \draw (1,0.95) -- (-0.2,0.95) node[left] {$\beta$}; 
      \draw[scale=0.5,domain=0.01:5.49,smooth,variable=\x,black,line width=1pt] plot ({\x},{0.60*\x}) node[right] {\small $F_1(\gamma)$}; 
       \draw[dashed] (1.051062,1.382) -- (1.051062,0) node[below] {$\gamma^\star$};
       \draw[scale=0.5,domain=0.01:5.49,smooth,variable=\x,black,line width=1pt] plot ({\x},{4*(1-exp(-\x))/\x}) node[right] {\small $F_2(\gamma)$};
       \draw[scale=0.5,domain=0.01:5.49,smooth,variable=\x,black,line width=1pt] plot ({\x},{4*(1-exp(-\x))/\x + 0.60*\x}) node[right] {\small $F(\gamma)$};
      %     \draw[scale=0.5,domain=-1.92:1,smooth,variable=\x,blue] plot ({\x},{1/3*1/(2+\x)});
  %    \end{axis}
    \end{tikzpicture}
    \caption{\label{F1}The functions  $F(\gamma)$, $F_1(\gamma)$, and $F_2(\gamma)$ for $2\kappa_{\rm B}\geqslant \lambda t \,\kappa_{\rm NC}$ (left panel) and for $2\kappa_{\rm B}< \lambda t \,\kappa_{\rm NC}$ (right panel). In the former case $\gamma^\star=0$, whereas in the latter case $\gamma^\star>0$. }
    \end{figure}}

We can also find the (reduced) second moments $v_{ij}(t)$, with $i,j=1,2$, by using the techniques from Section \ref{sec:SM}. Indeed, for $t\geqslant 0$, we have
\begin{align*}
    v_{11}'(t)&=2\lambda \,m_1(t) - \gamma v_{11}(t),\\
    v_{12}'(t)&=\lambda \,m_2(t) -\gamma v_{12}(t),\\
    v_{22}'(t)&= \gamma v_{11}(t) + 2 \gamma \, v_{12}(t).
\end{align*}
Solving the differential equations we obtain
\begin{align*}
    v_{11}(t) &= \Big(\frac{\lambda}{\gamma}\Big)^2 \Big[ 2\,(1-e^{-\gamma t}) -2 \, \gamma \, t \, e^{- \gamma t} \Big]\\
    v_{12}(t) &= \Big(\frac{\lambda}{\gamma}\Big)^2 \Big[ \gamma  t \,(1+ e^{-\gamma t}) + 2( e^{-\gamma t} -1) \Big]\\
    v_{22}(t) &= \Big(\frac{\lambda}{\gamma}\Big)^2 \Big[2 \,(1- e^{-\gamma t}) + \gamma^2 \, t^2 - 2 \, \gamma \, t  \Big].
\end{align*}
Using the first and (reduced) second moments, we can derive the following expressions for the variances:
\begin{align*}
    {\mathbb V}\mathrm{ar}\, M_1(t) &= \Big(\frac{\lambda}{\gamma}\Big)^2 \Big[-2 \, \gamma \, t \, e^{-\gamma t} - e^{-2 \gamma t} - \frac{\gamma}{\lambda}\, e^{-\gamma t} + 1 + \frac{\gamma}{\lambda} \Big]\\
    {\mathbb V}\mathrm{ar}\, M_2(t) &= \Big(\frac{\lambda}{\gamma}\Big)^2 \, \Big[-2 \, \gamma \, t \, e^{-\gamma t} - e^{-2 \gamma t} + \frac{\gamma}{\lambda}\, e^{-\gamma t} +\frac{\gamma^2}{\lambda} \, t + 1 - \frac{\gamma}{\lambda} \Big].
\end{align*}

\subsection{Faulty upload link}
This model can be further refined in various ways. In the first place, the link between the clients and the central storage may be faulty, in that it alternates between ``functioning'' and ``broken''. Suppose that the time the link is up (down, respectively) is exponentially distributed with parameter $q_{\rm U}$ ($q_{\rm D}$, respectively). Then it easily seen, using the techniques developed in Section \ref{sec:SM}, that, in self-evident notation, with the second index corresponding to the link being up or down,
\begin{align*}
    m_{1{\rm U}}'(t) &= -q_{\rm U}m_{1{\rm U}}(t) + q_{\rm D}m_{1{\rm D}}(t) +\lambda\pi_{\rm U}(t) - \gamma m_{1{\rm U}}(t),\\
     m_{1{\rm D}}'(t) &= q_{\rm U}m_{1{\rm U}}(t) - q_{\rm D}m_{1{\rm D}}(t) +\lambda\pi_{\rm D}(t),\\
      m_{2{\rm U}}'(t) &= -q_{\rm U}m_{2{\rm U}}(t) + q_{\rm D}m_{2{\rm D}}(t) + \gamma m_{1{\rm U}}(t),\\
       m_{2{\rm D}}'(t) &= q_{\rm U}m_{2{\rm U}}(t) - q_{\rm D}m_{2{\rm D}}(t);
\end{align*}
here the probabilities $\pi_{\rm U}(t)$ and $\pi_{\rm D}(t)$ follow from
\[\pi_{\rm U}(t) = \pi_{\rm U}(0)\pi_{\rm UU}(t) + \pi_{\rm D}(0)\pi_{\rm DU}(t),\:\:\:
\pi_{\rm D}(t) = \pi_{\rm U}(0)\pi_{\rm UD}(t) + \pi_{\rm D}(0)\pi_{\rm DD}(t),\]
with, abbreviating $q:=q_{\rm D}+q_{\rm U}$,
\[
   \pi_{\rm DU}(t)= 1- \pi_{\rm DD}(t) = 
\frac{q_{\rm D}}{q}(1-e^{-qt}),\:\:\:\:
  \pi_{\rm UD}(t)= 1- \pi_{\rm UU}(t) = 
\frac{q_{\rm U}}{q}(1-e^{-qt}).\]
The (reduced) second moments can be found relying on the theory of Section \ref{sec:SM}; we do not provide the differential equations here.

\subsection{Failures in storage unit}
A second extension includes the effect that the storage system itself can also be faulty, where upon failure all stored files are lost. We let these failure occur according to a Poisson process with intensity $\bar\gamma$. Let $M_3(t)$ denote the number of lost files. While the differential equations for $m_{1{\rm U}}(t)$ and $m_{1{\rm D}}(t)$ remain unchanged, the other ones become
\begin{align*}
    m_{2{\rm U}}'(t) &= -q_{\rm U}m_{2{\rm U}}(t) + q_{\rm D}m_{2{\rm D}}(t) + \gamma m_{1{\rm U}}(t) -\bar\gamma m_{2{\rm U}}(t)  ,\\
       m_{2{\rm D}}'(t) &= q_{\rm U}m_{2{\rm U}}(t) - q_{\rm D}m_{2{\rm D}}(t)-\bar\gamma m_{2{\rm D}}(t);\\
        m_{3{\rm U}}'(t) &= -q_{\rm U}m_{3{\rm U}}(t) + q_{\rm D}m_{3{\rm D}}(t)+\bar\gamma m_{2{\rm U}}(t), \\
        m_{3{\rm D}}'(t) &= q_{\rm U}m_{3{\rm U}}(t) - q_{\rm D}m_{3{\rm D}}(t)+\bar\gamma m_{2{\rm D}}(t).
\end{align*}
For a company offering this service, a typical design goal would be: how small should $q_{\rm U}$, respectively $\bar\gamma$, be to make sure that the mean number of lost files at time $t$, i.e., $m_{3{\rm U}}(t) +
        m_{3{\rm D}}(t)$, is below a given threshold? More advanced criteria could also involve the corresponding variance. 

To mitigate the effect of file loss due to failures of the storage unit,  an evident policy is to copy all files to {\it multiple} storage units. For instance, one could consider the mechanism in which every file is simultaneously copied to two storage units at Poisson instants (again with rate~$\gamma$), and that each of the storage units  fails at Poisson instants (again with rate~$\bar\gamma$). This means that a file is lost only if both storage units have failed. The above analysis (for the case of a single storage unit, that is) can be extended in an evident manner to this situation, thus facilitating the quantification of the performance gain due to making multiple backups. In order to decide whether such a policy should be implemented, this gain should be compared to the cost of the additional storage unit. 

\section{Concluding remarks} \label{Sec:CR}
This paper developed a versatile Markovian model that describes the dissemination of wealth over a population, but with the potential to be broadly applied across a wide range of disciplines. It is demonstrated that the evaluation of transient moments reduces to solving a system of coupled linear differential equations, while their stationary counterparts require solving (ordinary) linear systems of equations. Also a stability condition has been developed. Three examples evidence the model's wide application potential.

In various directions there is scope for generalizations. It is in particular noted that, conditional on the state of the background process, the wealth units evolve independently, which is not in all application areas a realistic mechanism. In addition, one may try to relax the various exponentiality assumptions. 

\section*{Declarations}

\subsection*{Conflicts of Interest} The authors declare that they have no conflicts of interest.

\subsection*{Availability of Data and Materials} Data: not applicable. Materials: the Python code used in generating the figures in this article are available from the corresponding author on reasonable request.

\bibliographystyle{plain}

\end{document}